\documentclass[preprint,authoryear,12pt]{elsarticle}

\usepackage{amssymb}
\def\th#1{\vspace{1mm}\noindent{\bf #1}}
\def\pf#1{\vspace{1mm}\noindent{\bf #1}}

\journal{Stochastic Processes and Their Applications}

\begin{document}

\begin{frontmatter}



\title{A central limit theorem under sublinear expectations\tnoteref{label1}}
\tnotetext[label1]{This work was partially supported by the
National Natural Science Foundation of China (Grant No. 10771122),
Natural Science Foundation of Shandong Province of China (Grant
No. Y2006A08) and National Basic Research Program of China (Grant
No. 2007CB814900).}

\author{Min Li}

\author{Yufeng Shi\corref{cor1}}
\ead{yfshi@sdu.edu.cn.}

\cortext[cor1]{corresponding author.}

\address{School of Mathematics, Shandong University,
Jinan 250100, China}

\begin{abstract}
In this paper we consider a sequence of random variables with mean
uncertainty in a sublinear expectation space. Without the
hypothesis of identical distributions, we show a new central limit
theorem under the sublinear expectations.
\end{abstract}


\end{frontmatter}


\section{Introduction}
Motivated by the coherent risk measures (cf. [1, 2]) and uncertain
volatility models in finance (see, e.g. [11]), Peng [12, 14]
recently introduced the notion of sublinear expectation which is
not based on a classical probability space. Under the sublinear
expectations, a random variable $X$ in a sublinear expectation
space $\left( \Omega ,{\cal H},\hat E\right) $ (for its
definition, see Section 2 of this paper) is said to be of
$G$-normal distribution with zero mean (cf. [14, 17]), if for each
$Y$ which is an independent copy of $X$, it holds that
\[
aX+bY\stackrel{d}{=}\sqrt{a^2+b^2}X,\ \quad \forall a,b\geq 0.
\]
Just as the classical normal distributions in probability theory, for a $G$%
-normal distributed random variable $X$, we have (cf. [14])
\[
\hat E\left[ \varphi \left( X\right) \right] =u\left( 1,0\right)
,\quad \forall \varphi \in C_{b,Lip}(R),
\]
where $u\left( t,x\right) $ is the unique viscosity solution for
the following heat equation
\[
\left\{
\begin{array}{c}
\partial _tu-G\left( \partial _{xx}^2u\right) =0,\\
u\mid_{t=0}=\varphi ,
\end{array}
\right.
\]
where $G\left( \alpha \right) :=\hat E\left[ \frac 12\alpha
X^2\right] $. In the theory of sublinear expectations, the above
heat equation often plays a role of characteristic function in
probability theory.

On the basis of $G$-normal distribution, $G$-Brownian motion can
be defined, and the corresponding stochastic calculus with respect
to the $G$-Brownian motions and the related It\^o's formula can
also be established (cf. [12, 14]). Since the importance of law of
large numbers (LLN) and central limit theorem (CLT) in probability
theory, Peng [13, 15] has shown the
corresponding LLN and CLT under sublinear expectations, which indicate that $%
G$-normal distributions play the same important role in the theory
of sublinear expectations as the normal distributions in the
classical probability theory.

Due to the significance of sublinear expectations in finance and
statistics, the theory of sublinear expectations has been
attracting more and more attentions in both pure and applied
mathematics (see, e.g. [6], [8], [18], [19] and [20]).

The purpose of this paper is to investigate one of the very
important fundamental results in the theory of sublinear
expectations---Central Limit Theorem. Until now, all the results
on central limit theorems under sublinear expectations require
that the sequence of random variables is independent and
identically distributed. Analogous to the CLT in the probability
theory, a natural question is whether one can weaken the
hypothesis of identical distributions for the CLT under sublinear
expectations?

In this paper, without the hypothesis of identical distributions,
we prove a new central limit theorem within mean uncertainty under
the sublinear expectations, which extends Peng's results.
Precisely, for a sequence of random variables $\left\{ \left(
X_n,Y_n\right) \right\} _{n=1}^\infty $ in a sublinear expectation
space, we only require the random variable $\left(
X_{n+1},Y_{n+1}\right) $ is independent to $\left\{ \left(
X_1,Y_1\right) ,\left( X_2,Y_2\right) ,\cdot \cdot \cdot ,\left(
X_n,Y_n\right) \right\} $, $n=1,2,3,\cdots $.

\section{Basic settings}
Let $\Omega $ be a given set and ${\cal H}$ a vector lattice
of real functions defined on $\Omega $, including $1$, such that if $%
X_1,\cdot \cdot \cdot ,X_n\in {\cal H}$ then $\varphi \left(
X_1,\cdot \cdot \cdot ,X_n\right) \in {\cal H}$ for each $\varphi
\in C_{b,Lip}\left( R^n\right) $ where $C_{b,Lip}\left( R^n\right)
$ denotes the space of bounded and Lipschitz continuous functions
defined on $R^n$. ${\cal H}$ is considered as a space of ``random
variables''. We denote by $\left\langle x,y\right\rangle $ the
scalar product of $x,y\in R^n$ and by $\left| x\right|
=\left\langle x,x\right\rangle ^{1/2}$ the Euclidean norm of $x\in
R^n$. Let $S\left( n\right) $ be the collection of $n\times
n$-symmetric matrices. $S\left( n\right) $ is obviously a Hilbert
space with the scalar product $\left\langle P,Q\right\rangle
=tr\left[ PQ\right] $.

Now we give some related definitions about sublinear expectations
(see [12-17] for details).

\th {Definition 2.1.}\ A sublinear expectation $\hat E$ on ${\cal H}$ is a functional $\hat E:%
{\cal H}\rightarrow R$ satisfying the following properties: for
all $X,Y\in {\cal H}$, we have
\begin{enumerate} [(a)]
\item Monotonicity: If $X\geq Y$, then $\hat E\left[ X\right] \geq
\hat E\left[ Y\right] .$
\item Constant preserving: $\hat E\left[ c\right] =c$, $\forall
c\in R.$
\item Sub-additivity: $\hat E\left[ X+Y\right] \leq \hat
E\left[ X\right] +\hat E\left[ Y\right] .$
\item Positive homogeneity: $\hat E\left[ \lambda X\right]
=\lambda \hat E\left[ X\right] ,\forall \lambda \geq 0.$
\end{enumerate}
The triple $\left( \Omega ,{\cal H},\hat E\right) $ is called a  %
sublinear expectation space.

\th {Definition 2.2.}\ Let $X_1$ and $X_2$ be two $n$-dimensional
random vectors defined in sublinear expectation spaces $\left(
\Omega _1,{\cal H}_1,\hat E_1\right) $ and $\left( \Omega _2,{\cal
H}_2,\hat E_2\right) $, respectively. They are called identically
distributed, denoted by $X_1\stackrel{d}{=}X_2$, if
\[
\hat E_1\left[ \varphi \left( X_1\right) \right] =\hat E_2\left[
\varphi \left( X_2\right) \right] ,\quad \forall \varphi \in
C_{b,Lip}\left( R^n\right) .
\]

\th {Definition 2.3.}\ In a sublinear expectation space $\left(
\Omega ,{\cal H},\hat E\right) $ a random vector $Y=\left(
Y_1,\cdot \cdot \cdot ,Y_n\right) $, $Y_i\in {\cal H} $ is said to
be independent to another random vector $X=\left( X_1,\cdot \cdot
\cdot ,X_m\right) $, $X_i\in {\cal H}$ under $\hat E\left[ \cdot
\right] $ if for each test function $\varphi \in C_{b,Lip}\left(
R^m\times R^n\right) $ we have
\[
\hat E\left[ \varphi \left( X,Y\right) \right] =\hat E\left[ \hat
E\left[ \varphi \left( x,Y\right) \right] _{x=X}\right] .
\]

The definition of independence means that any realization of $X$
does not change the distributional uncertainty of $Y$. At the same
time, the fact
that $Y$ is independent to $X$ does not imply that $X$ is independent to $Y$%
, an example can be seen in [14] or [15].

\th {Definition 2.4.}\ ($G$-normal distribution with mean
uncertainty) A pair of $d$-dimensional random vectors $\left(
X,Y\right) $ in a sublinear expectation space $\left( \Omega
,{\cal H},\hat E\right) $ is called $G$-normal distributed if for
each $a,b\geq 0$ we have
\[
\left( aX+b\overline{X},a^2Y+b^2\overline{Y}\right)
\stackrel{d}{=}\left( \sqrt{a^2+b^2}X,\left( a^2+b^2\right)
Y\right) ,\quad \forall a,b\geq 0,
\]
where $\left( \overline{X},\overline{Y}\right) $ is an independent copy of $%
\left( X,Y\right) $.

The following propositions (cf. [12, 15]) play an important role
in this paper.

\th {Proposition 2.5.}\ {\it Let $G:R^d\times S\left( d\right)
\rightarrow R$ be a given sublinear functional continuous in
$\left( 0,0\right) $ and satisfying the following properties: for
any $\left( p,A\right) $ and $\left( \bar p,\bar A\right) \in
R^d\times S\left( d\right) $
\[
\left\{
\begin{array}{c}
G\left( p+\overline{p},A+\overline{A}\right) \leq G\left(
p,A\right)
+G\left( \overline{p},\overline{A}\right) , \\
G\left( \lambda p,\lambda A\right) =\lambda G\left( p,A\right)
,\quad
\forall \lambda \geq 0, \\
G\left( p,A\right) \geq G\left( p,\overline{A}\right) ,\quad if
A\geq \overline{A}.
\end{array}
\right.
\]
Then there exists a pair of $d$-dimensional $G$-normal distributed
random vectors $\left( X,Y\right) $ in some sublinear expectation
space $\left( \Omega ,{\cal H},\hat E\right) $ such that
\[
G\left( p,A\right) =\hat E\left[ \frac 12\left\langle
AX,X\right\rangle +\left\langle p,Y\right\rangle \right] ,\quad
\forall \left( p,A\right) \in R^d\times S\left( d\right) .
\]}

\th {Proposition 2.6.}\ {\it Let $\left( X,Y\right) $ be a
$G$-normal distributed random vector in a sublinear expectation
space $\left( \Omega ,{\cal H},\hat E\right) $. For each $\varphi
\in C_{b,Lip}\left( R^d\right) $ we define a function
\[
v\left( t,x\right) :=\hat E\left[ \varphi \left(
x+\sqrt{t}X+tY\right) \right] ,\quad \forall \left( t,x\right) \in
\left[ 0,\infty \right) \times R^d.
\]
Then $v$ is the unique viscosity solution of the following
parabolic partial differential equation (PDE)
\[
\partial _tv-G\left( D_xv,D_x^2v\right) =0,\quad v\mid _{t=0}=\varphi ,
\]
where
\[
G\left( p,A\right) =\hat E\left[ \frac 12\left\langle
AX,X\right\rangle +\left\langle p,Y\right\rangle \right] ,\quad
\forall \left( p,A\right) \in R^d\times S\left( d\right) .
\]}

\th {Proposition 2.7.}\ {\it Let $X,Y$ be two random variables in
a sublinear expectation space $\left( \Omega ,{\cal H},\hat
E\right) $, then for $1<p,q<\infty ,\frac 1p+\frac 1q=1 $, we have
\[
\hat E\left[ \left| XY\right| \right] \leq \left( \hat E\left[
\left| X\right| ^p\right] \right) ^{1/p}\cdot \left( \hat E\left[
\left| Y\right| ^q\right] \right) ^{1/q},
\]
In particular, for $1\leq p\leq p^{^{\prime }}$, we have
\[
\left( \hat E\left[ \left| X\right| ^p\right] \right) ^{1/p}\leq
\left( \hat E\left[ \left| X\right| ^{p^{^{\prime }}}\right]
\right) ^{1/p^{^{\prime }}}.
\]}

\section{Main results}
Now we give the main result in this paper---Central Limit Theorem.
For the simplicity of notations, we first prove the
$1$-dimensional case of CLT.

\th {Theorem 3.1.}\ {\it In a sublinear expectation space $\left(
\Omega ,{\cal H},\hat E\right) $, let $\left\{ \left(
X_i,Y_i\right) \right\} _{i=1}^\infty $ be a sequence of $R\times
R$-valued random variables and $\left( \xi ,\zeta \right) $ be a
pair of $G$-normal distributed random variables. We assume that
\begin{enumerate} [(i)]
\item $\left( X_{i+1},Y_{i+1}\right) $ is independent to $\left\{
\left(
X_1,Y_1\right) ,\cdot \cdot \cdot ,\left( X_i,Y_i\right) \right\} ,$ for $%
i=1,2,\cdot \cdot \cdot $;
\item $\hat E\left[ X_i\right] =\hat E\left[ -X_i\right] =0,$ $\hat
E\left[ \left| X_i\right| ^3\right] \leq M,$ $\hat E\left[ \left|
Y_i\right| ^3\right] \leq M,$ where $M$ is a positive constant;
\item $\lim\limits_{n\rightarrow \infty }\frac 1n\sum\limits_{i=1}^{n}
\hat E\left[ \left| X_i^2-\xi ^2\right| ^2\right] =0, $$
\lim\limits_{n\rightarrow \infty }\frac
1n\sum\limits_{i=1}^{n}\hat E\left[ \left| Y_i-\zeta \right|
^2\right] =0;$
\item there exists $\beta >0,$ such that $\hat E\left[ a\xi
^2\right] -\hat E\left[ \overline{a}\xi ^2\right] \geq \beta
\left( a-\overline{a}\right) ,$ for any $a,\overline{a}\in R$ with
$a\geq \overline{a}.$
\end{enumerate}
Then the sequence $\left\{
\frac{S_n}{\sqrt{n}}+\frac{T_n}n\right\} _{n=1}^\infty $, where
$S_n=X_1+\cdot \cdot \cdot +X_n,$ $T_n=Y_1+\cdot \cdot \cdot
+Y_n,$ converges in law to $\xi +\zeta $:
\[
\lim\limits_{n\rightarrow \infty }\hat E\left[ \varphi \left( \frac{S_n%
}{\sqrt{n}}+\frac{T_n}n\right) \right] =\hat E\left[ \varphi
\left( \xi +\zeta \right) \right] ,\quad \forall \varphi \in
C_{b,Lip}\left( R\right) ,
\]
where the sublinear function $G:R\times R\rightarrow R$ is defined
by
\[
G\left( p,a\right) :=\hat E\left[ p\zeta +\frac 12a\xi ^2\right] .
\]}

\pf{Proof.}\ For any $\varphi \in C_{b,Lip}\left( R\right) $, and
a small but fixed $h>0$, let $V$ be the unique viscosity solution
of
\begin{equation}
\partial _tV+G\left( \partial _xV,\partial _{xx}^2V\right) =0,\quad \left(
t,x\right) \in \left[ 0,1+h\right] \times R,\quad V\mid
_{t=1+h}=\varphi .
\end{equation}
Since $\left( \xi ,\zeta \right) $ is of $G$-normal distribution,
from Proposition 6, we have
\[
V\left( t,x\right) =\hat E\left[ \varphi \left( x+\sqrt{1+h-t}\xi
+\left( 1+h-t\right) \zeta \right) \right] .
\]
Particularly,
\begin{equation}
V\left( h,0\right) =\hat E\left[ \varphi \left( \xi +\zeta \right)
\right] ,\quad V\left( 1+h,x\right) =\varphi \left( x\right) .
\end{equation}
Since (1) is a uniformly parabolic PDE and $G$ is a convex
function, thus, by the interior regularity of $V$ (see Krylov
[10], Theorem 6.2.3), we have
\[
\left| V\right| _{C^{1+\alpha /2,2+\alpha }\left( \left[
0,1\right] \times R\right) }<\infty ,
\]
for some $\alpha \in \left(0,1\right) $.\\
We set $\delta =\frac 1n$, $S_0=T_0=0$ and $\overline{S}_i=\sqrt{\delta }%
S_i+\delta T_i$, then
\begin{eqnarray*}
&&V\left( 1,\sqrt{\delta }S_n+\delta T_n\right) -V\left( 0,0\right) \\
&=&V\left( 1,\overline{S}_n\right) -V\left( 0,0\right) \\
&=&\sum\limits_{i=0}^{n-1}\left\{ V\left( \left( i+1\right) \delta
,\overline{S}_{i+1}\right) -V\left( i\delta ,\overline{S}_i\right)
\right\} \\
&=&\sum\limits_{i=0}^{n-1}\left\{ \left[ V\left( \left(
i+1\right) \delta ,\overline{S}_{i+1}\right) -V\left( i\delta ,\overline{S}%
_{i+1}\right) \right] +\left[ V\left( i\delta
,\overline{S}_{i+1}\right)
-V\left( i\delta ,\overline{S}_i\right) \right] \right\} \\
&=&\sum\limits_{i=0}^{n-1}\left\{ I_\delta ^i+J_\delta ^i\right\}
,
\end{eqnarray*}
with, by Taylor's expansion,
\begin{eqnarray*}
J_\delta ^i &=&\partial _tV\left( i\delta ,\overline{S}_i\right)
\delta +\frac 12\partial _{xx}^2V\left( i\delta
,\overline{S}_i\right) X_{i+1}^2\delta +\partial _xV\left( i\delta
,\overline{S}_i\right) \left(
X_{i+1}\sqrt{\delta }+Y_{i+1}\delta \right) , \\
I_\delta ^i &=&\int_0^1\left( \partial _tV\left( i\delta +\beta \delta ,%
\overline{S}_{i+1}\right) -\partial _tV\left( i\delta ,\overline{S}%
_{i+1}\right) \right) d\beta \delta +\left( \partial _tV\left( i\delta ,%
\overline{S}_{i+1}\right) -\partial _tV\left( i\delta
,\overline{S}_i\right)
\right) \delta \\
&&\ \ +\int_0^1\int_0^1\left[ \partial _{xx}^2V\left( i\delta ,\overline{S}%
_i+\gamma \beta \left( X_{i+1}\sqrt{\delta }+Y_{i+1}\delta \right)
\right) -\partial _{xx}^2V\left( i\delta ,\overline{S}_i\right)
\right] \beta d\gamma d\beta \\
&&\ \times\left( X_{i+1}\sqrt{\delta }+Y_{i+1}\delta \right) ^2 \
+\frac 12\partial _{xx}^2V\left( i\delta ,\overline{S}_i\right)
\left( Y_{i+1}^2\delta ^2+2X_{i+1}Y_{i+1}\delta ^{3/2}\right) .
\end{eqnarray*}
Thus
\begin{eqnarray}
\hat E\left[ \sum\limits_{i=0}^{n-1}J_\delta ^i\right] -\hat
E\left[ -\sum\limits_{i=0}^{n-1}I_\delta ^i\right] &\leq &\hat
E\left[ V\left( 1,\overline{S}_n\right) \right] -V\left(
0,0\right)
\nonumber \\
&\leq &\hat E\left[ \sum\limits_{i=0}^{n-1}J_\delta ^i\right]
+\hat E\left[ \sum\limits_{i=0}^{n-1}I_\delta ^i\right] .
\end{eqnarray}
For the 3rd term of $J_\delta ^i$, by (i) and (ii) we have
\[
\widehat{E}\left[ \partial _xV\left( i\delta ,\overline{S}_i\right) X_{i+1}%
\sqrt{\delta }\right] =\widehat{E}\left[ -\partial _xV\left( i\delta ,%
\overline{S}_i\right) X_{i+1}\sqrt{\delta }\right] =0.
\]
We then combine the above equality with (1) as well as the
independence of $\left( X_{i+1},Y_{i+1}\right) $ to $\left\{
\left( X_1,Y_1\right) ,\cdot \cdot \cdot ,\left( X_i,Y_i\right)
\right\} $, it follows that
\begin{eqnarray*}
\hat E\left[ J_\delta ^i\right] &=&\hat E\left[ \partial _tV\left( i\delta ,%
\overline{S}_i\right) \delta +\frac 12\partial _{xx}^2V\left( i\delta ,%
\overline{S}_i\right) X_{i+1}^2\delta +\partial _xV\left( i\delta ,\overline{%
S}_i\right) Y_{i+1}\delta \right] \\
\ &=&\hat E\left[ \partial _tV\left( i\delta
,\overline{S}_i\right) \delta +\hat E\left[ \frac 12\partial
_{xx}^2V\left( i\delta ,\overline{S}_i\right) X_{i+1}^2\delta
+\partial _xV\left( i\delta ,\overline{S}_i\right)
Y_{i+1}\delta \right] \right] \\
\ &\leq &\delta \hat E\left[ \partial _tV\left( i\delta ,\overline{S}%
_i\right) +\hat E\left[ \frac 12\partial _{xx}^2V\left( i\delta ,\overline{S}%
_i\right) \xi ^2+\partial _xV\left( i\delta ,\overline{S}_i\right)
\zeta
\right] \right] \\
&&\ \ +\delta \hat E\left[ \frac 12\partial _{xx}^2V\left( i\delta ,%
\overline{S}_i\right) \left( X_{i+1}^2-\xi ^2\right) +\partial
_xV\left(
i\delta ,\overline{S}_i\right) \left( Y_{i+1}-\zeta \right) \right] \\
\ &=&\delta \hat E\left[ \frac 12\partial _{xx}^2V\left( i\delta ,\overline{S%
}_i\right) \left( X_{i+1}^2-\xi ^2\right) +\partial _xV\left( i\delta ,%
\overline{S}_i\right) \left( Y_{i+1}-\zeta \right) \right] \\
\ &\leq &\frac \delta 2\hat E\left[ \partial _{xx}^2V\left( i\delta ,%
\overline{S}_i\right) \left( X_{i+1}^2-\xi ^2\right) \right]
+\delta \hat E\left[ \partial _xV\left( i\delta
,\overline{S}_i\right) \left(
Y_{i+1}-\zeta \right) \right] \\
\ &\leq &\frac \delta 2\left( \hat E\left[ \left| \partial
_{xx}^2V\left( i\delta ,\overline{S}_i\right) \right| ^2\right]
\right) ^{1/2}\cdot \left(
\hat E\left[ \left| X_{i+1}^2-\xi ^2\right| ^2\right] \right) ^{1/2} \\
&&\ \ +\delta \left( \hat E\left[ \left| \partial _xV\left( i\delta ,%
\overline{S}_i\right) \right| ^2\right] \right) ^{1/2}\cdot \left(
\hat E\left[ \left| Y_{i+1}-\zeta \right| ^2\right] \right)
^{1/2}.
\end{eqnarray*}
But since $\partial _{xx}^2V$ is uniformly $\frac \alpha
2$-H\"older continuous in $t$ and $\alpha $-H\"older continuous in
$x$ on $\left[ 0,1\right] \times R$, it follows that
\[
\left| \partial _{xx}^2V\left( i\delta ,\overline{S}_i\right)
-\partial _{xx}^2V\left( 0,0\right) \right| \leq C\left( \left|
\overline{S}_i\right| ^\alpha +\left| i\delta \right| ^{\frac
\alpha 2}\right) ,
\]
where $C$ is some positive constant. We note that
\begin{eqnarray*}
\hat E\left[ \left| \overline{S}_i\right| ^{2\alpha }\right] &\leq
&\hat
E\left[ \left| \overline{S}_i\right| ^{2\alpha }\vee 1\right] \\
\ &=&\hat E\left[ \left( \left| \overline{S}_i\right| \vee
1\right)
^{2\alpha }\right] \\
\ &\leq &\hat E\left[ \left( \left| \overline{S}_i\right| \vee
1\right)
^2\right] \\
\ &\leq &\hat E\left[ \left( \left| \overline{S}_i\right|
+1\right) ^2\right]
\\
\ &\leq &2\hat E\left[ \left| \overline{S}_i\right| ^2\right] +2 \\
\ &\leq &\frac 4n\hat E\left[ \left( \sum\limits_{j=1}^{i}%
X_j\right) ^2\right] +\frac 4{n^2}\hat E\left[ \left( \sum\limits_{j=1}^{i}%
Y_j\right) ^2\right] +2 \\
\ &\leq &\frac 4n\sum\limits_{j=1}^{n}\hat E\left[ X_j^2\right]
+\frac 4n\sum\limits_{j=1}^{n}\hat E\left[ Y_j^2\right] +2,
\end{eqnarray*}
at the same time,
\begin{eqnarray*}
\frac 1n\sum\limits_{j=1}^{n}\hat E\left[ Y_j^2\right] &=&\frac
1n\sum\limits_{j=1}^{n}\hat E\left[ \left(
Y_j-\zeta +\zeta \right) ^2\right] \\
\ &\leq &\frac 2n\sum\limits_{j=1}^{n}\hat E\left[ \left(
Y_j-\zeta \right) ^2\right] +2\hat E\left[ \zeta ^2\right] .
\end{eqnarray*}
From (iii), it follows that
\begin{eqnarray}
\lim\limits_{n\rightarrow \infty }\frac 1n\sum\limits_{i=1}^{n}%
\left( \hat E\left[ \left| X_i^2-\xi ^2\right| ^2\right] \right)
^{1/2} &=&0,   \\
\lim\limits_{n\rightarrow \infty }\frac 1n\sum\limits_{i=1}^{n}%
\left( \hat E\left[ \left| Y_{i}-\zeta \right| ^2\right] \right)
^{1/2} &=&0,   \\
\lim\limits_{n\rightarrow \infty }\frac 1n\sum\limits_{j=1}^{n}%
\hat E\left[ X_j^2\right] &=&\hat E\left[ \xi ^2\right] .
\nonumber
\end{eqnarray}
Then there exists a constant $C_1>0$ such that
\[
\left( \hat E\left[ \left| \partial _{xx}^2V\left( i\delta ,\overline{S}%
_i\right) \right| ^2\right] \right) ^{1/2}\leq 2C_1.
\]
Since $\varphi \in C_{b,Lip}\left( R\right) $, there exists a constant $%
C_2>0 $ such that
\[
\left| \varphi \left( x\right) -\varphi \left( y\right) \right|
\leq C_2\left| x-y\right| ,\quad \forall x,y\in R.
\]
Then $\forall t\in \left[ 0,1\right] $, $\forall x,y\in R$, we
have
\begin{eqnarray*}
&&\left| V\left( t,x\right) -V\left( t,y\right) \right| \\
&=&\left| \hat E\left[ \varphi \left( x+\sqrt{1+h-t}\xi +\left(
1+h-t\right) \zeta \right) \right] -\hat E\left[ \varphi \left(
y+\sqrt{1+h-t}\xi +\left(
1+h-t\right) \zeta \right) \right] \right| \\
&\leq &\hat E\left[ \left| \varphi \left( x+\sqrt{1+h-t}\xi
+\left( 1+h-t\right) \zeta \right) -\varphi \left(
y+\sqrt{1+h-t}\xi +\left(
1+h-t\right) \zeta \right) \right| \right] \\
&\leq &C_2\left| x-y\right| .
\end{eqnarray*}
From the above inequality, we obtain
\[
\left| \partial _xV\left( t,x\right) \right| \leq C_2,\quad
\forall \left( t,x\right) \in \left[ 0,1\right] \times R.
\]
Thus
\[
\hat E\left[ J_\delta ^i\right] \leq \frac{C_1}n\left( \hat
E\left[ \left| X_{i+1}^2-\xi ^2\right| ^2\right] \right)
^{1/2}+\frac{C_2}n\left( \hat E\left[ \left| Y_{i+1}-\zeta \right|
^2\right] \right) ^{1/2},
\]
it follows that
\begin{eqnarray*}
\hat E\left[ \sum\limits_{i=1}^{n-1}J_\delta ^i\right] &\leq
&\sum\limits_{i=1}^{n-1}\hat E\left[ J_\delta ^i\right] \\
\ &\leq &\frac{C_1}n\sum\limits_{i=1}^{n-1}\left( \hat
E\left[ \left| X_{i+1}^2-\xi ^2\right| ^2\right] \right) ^{1/2}+\frac{C_2}n%
\sum\limits_{i=1}^{n-1}\left( \hat E\left[ \left| Y_{i+1}-\zeta
\right| ^2\right] \right) ^{1/2},
\end{eqnarray*}
then by (4) and (5), we have
\[
\lim\limits_{n\rightarrow \infty }\hat E\left[ \sum\limits_{i=0}^{n-1}%
J_\delta ^i\right] \leq 0.
\]
Similarly, we also have
\begin{eqnarray*}
\hat E\left[ J_\delta ^i\right] &\geq &-\frac \delta 2\left( \hat
E\left[ \left| \partial _{xx}^2V\left( i\delta
,\overline{S}_i\right) \right| ^2\right] \right) ^{1/2}\left( \hat
E\left[ \left| X_{i+1}^2-\xi ^2\right|
^2\right] \right) ^{1/2} \\
&&\ \ -\delta \left( \hat E\left[ \left| \partial _xV\left( i\delta ,%
\overline{S}_i\right) \right| ^2\right] \right) ^{1/2}\left( \hat
E\left[ \left| Y_{i+1}-\zeta \right| ^2\right] \right) ^{1/2},
\end{eqnarray*}
and
\[
\lim\limits_{n\rightarrow \infty }\hat E\left[ \sum\limits_{i=0}^{n-1}%
J_\delta ^i\right] \geq 0.
\]
Thus
\begin{equation}
\lim\limits_{n\rightarrow \infty }\hat E\left[ \sum\limits_{i=0}^{n-1}%
J_\delta ^i\right] =0.
\end{equation}
For $I_\delta ^i$, since both $\partial _tV$ and $\partial
_{xx}^2V$ are uniformly $\frac \alpha 2$-H\"older continuous in
$t$ and $\alpha $-H\"older continuous in $x$ on $\left[ 0,1\right]
\times R$, we have
\begin{eqnarray*}
\left| I_\delta ^i\right| &\leq &C_3\delta ^{1+\alpha /2}\left(
1+\left|
X_{i+1}+\sqrt{\delta }Y_{i+1}\right| ^\alpha +\left| X_{i+1}+\sqrt{\delta }%
Y_{i+1}\right| ^{2+\alpha }\right) \\
&&\ \ +\frac 12\left| \partial _{xx}^2V\left( i\delta
,\overline{S}_i\right) \left( Y_{i+1}^2\delta
^2+2X_{i+1}Y_{i+1}\delta ^{3/2}\right) \right| ,
\end{eqnarray*}
where $C_3$ is some positive constant. From (i), it follows that
\begin{eqnarray*}
\hat E\left[ \left| \partial _{xx}^2V\left( i\delta
,\overline{S}_i\right) X_{i+1}Y_{i+1}\delta ^{3/2}\right| \right]
&=&\frac 1{n^{3/2}}\hat E\left[ \left| \partial _{xx}^2V\left(
i\delta ,\overline{S}_i\right) \right|
\right] \hat E\left[ \left| X_{i+1}Y_{i+1}\right| \right] \\
&\leq &\frac{2C_1}{n^{3/2}}\left( \hat E\left[ \left|
Y_{i+1}\right| ^3\right] \right) ^{1/3}\left( \hat E\left[ \left|
X_{i+1}\right|
^{3/2}\right] \right) ^{2/3} \\
&\leq &\frac{2C_1}{n^{3/2}}M^{1/3}\left( \hat E\left[ \left|
X_{i+1}\right|
^3\right] \right) ^{1/3} \\
&\leq &\frac{2C_1}{n^{3/2}}M^{2/3},
\end{eqnarray*}
and
\begin{eqnarray*}
\hat E\left[ \frac 12\left| \partial _{xx}^2V\left( i\delta ,\overline{S}%
_i\right) Y_{i+1}^2\delta ^2\right| \right] \leq
\frac{C_1}{n^2}M^{2/3}.
\end{eqnarray*}
At the same time, we can claim that
\begin{eqnarray*}
\hat E\left[ \left| X_{i+1}+\sqrt{\delta }Y_{i+1}\right|
^{2+\alpha }\right]
\leq \hat E\left[ \left| X_{i+1}+\sqrt{\delta }Y_{i+1}\right| ^3\right] ^{%
\frac{2+\alpha }3}\leq \left( 8M\right) ^{\frac{2+\alpha }3}.
\end{eqnarray*}
Therefore it follows that
\[
\lim\limits_{n\rightarrow \infty }\hat E\left[ \sum\limits_{i=0}^{n-1}%
I_\delta ^i\right] \leq \lim\limits_{n\rightarrow \infty }
\sum\limits_{i=0}^{n-1}\hat E\left[ \left| I_\delta ^i\right|
\right] =0.
\]
It also holds that
\[
\lim\limits_{n\rightarrow \infty }-\hat E\left[ -\sum\limits_{i=0}^{n-1}%
I_\delta ^i\right] \geq \lim\limits_{n\rightarrow \infty }
-\sum\limits_{i=0}^{n-1}\hat E\left[ \left| I_\delta ^i\right|
\right] =0.
\]
Thus
\begin{equation}
\lim\limits_{n\rightarrow \infty }\hat E\left[ \sum\limits_{i=0}^{n-1}%
I_\delta ^i\right] =\lim\limits_{n\rightarrow \infty }-\hat
E\left[ -\sum\limits_{i=0}^{n-1}I_\delta ^i\right] =0.
\end{equation}
Therefore from (3), (6) and (7), we have
\[
\lim\limits_{n\rightarrow \infty }\hat E\left[ V\left( 1,\sqrt{\delta }%
S_n+\delta T_n\right) \right] =V\left( 0,0\right) .
\]
On the other hand, for each $t,t^{^{\prime }}\in \left[ 0,1+h\right] $ and $%
x\in R$, we have
\[
\left| V\left( t,x\right) -V\left( t^{^{\prime }},x\right) \right|
\leq C\left( \sqrt{\left| t-t^{^{\prime }}\right| }+\left|
t-t^{^{\prime }}\right| \right) .
\]
Thus $\left| V\left( h,0\right) -V\left( 0,0\right) \right| \leq
C\left( \sqrt{h}+h\right) $ and
\begin{eqnarray*}
\left| \hat E\left[ V\left( 1,\overline{S}_n\right) \right] -\hat
E\left[ \varphi \left( \overline{S}_n\right) \right] \right|
=\left| \hat E\left[
V\left( 1,\overline{S}_n\right) \right] -\hat E\left[ V\left( 1+h,\overline{S%
}_n\right) \right] \right| \leq C\left( \sqrt{h}+h\right) .
\end{eqnarray*}
Then we can claim that
\[
\overline{\lim\limits_{n\rightarrow \infty }}\left| \hat E\left[
\varphi \left( \sqrt{\delta }S_n+\delta T_n\right) \right] -\hat
E\left[
\varphi \left( \xi +\zeta \right) \right] \right| \leq 2C\left( \sqrt{h}%
+h\right) .
\]
Since $h$ can be arbitrarily small, we obtain that
\[
\lim\limits_{n\rightarrow \infty }\hat E\left[ \varphi \left( \sqrt{%
\delta }S_n+\delta T_n\right) \right] =\hat E\left[ \varphi \left(
\xi +\zeta \right) \right] .
\]
The proof is completed.

It is not difficult to obtain the following statement.

\th {Remark 3.2.}\ In addition of the assumptions of Theorem 3.1, if there is a sequence of $%
R\times R$-valued random variables $\left\{ \left( \overline{X}_i,\overline{Y%
}_i\right) \right\} _{i=1}^\infty $ in another sublinear expectation space $%
\left( \Omega _1,{\cal H}_1,\widetilde{E}\right) $, such that
$\left( \overline{X}_{i+1},\overline{Y}_{i+1}\right) $ is
independent to $\left\{ \left(
\overline{X}_1,\overline{Y}_1\right) ,\cdot \cdot \cdot ,\left(
\overline{X}_i,\overline{Y}_i\right) \right\} $ and $\left( \overline{X}_i,%
\overline{Y}_i\right) \stackrel{d}{=}\left( X_i,Y_i\right) $, for $%
i=1,2,\cdot \cdot \cdot $, then we also have
\[
\lim\limits_{n\rightarrow \infty }\widetilde{E}\left[ \varphi
\left(
\frac{\sum\limits_{i=1}^{n}\overline{X}_i}{\sqrt{n}}+\frac{%
\sum\limits_{i=1}^{n}\overline{Y}_i}n\right) \right] =\hat E\left[
\varphi \left( \xi +\zeta \right) \right] ,\quad \forall \varphi
\in C_{b,Lip}\left( R\right) .
\]

By the same arguments, we can claim the multi-dimensional case of
CLT.

\th {Theorem 3.3.}\ {\it In a sublinear expectation space
$\left( \Omega ,{\cal H},\hat E\right) $, let $\left\{ \left(
X_i,Y_i\right) \right\} _{i=1}^\infty $ be a sequence of
$R^d\times R^d$-valued random vectors and $\left( \xi ,\zeta
\right) $ be a pair of $G$-normal distributed $d$-dimensional
random vectors. We assume that
\begin{enumerate} [(i)]
\item $\left( X_{i+1},Y_{i+1}\right) $ is independent to $\left\{
\left(
X_1,Y_1\right) ,\cdot \cdot \cdot ,\left( X_i,Y_i\right) \right\} ,$ for $%
i=1,2,\cdot \cdot \cdot $;
\item $\hat E\left[ X_i\right] =\hat E\left[ -X_i\right] =0$ and
$\hat E\left[ \left| X_i\right| ^3\right] \leq M,$ $\hat E\left[
\left| Y_i\right| ^3\right] \leq M,$ where $M$ is a positive
constant;
\item $\lim\limits_{n\rightarrow \infty }\frac 1n\sum\limits_{i=1}^{n}%
\hat E\left[ \left| X_iX_i^T-\xi \xi ^T\right|
^2\right] =0$ and $\lim\limits_{n\rightarrow \infty }\frac 1n\sum\limits_{i=1}^{n}%
\hat E\left[ \left| Y_i-\zeta \right| ^2\right] =0;$
\item there exists $\beta >0,$ such that $\hat E\left[ \left\langle
A\xi ,\xi \right\rangle \right] -\hat E\left[ \left\langle
\overline{A}\xi ,\xi \right\rangle \right] \geq \beta tr\left[
A-\overline{A}\right] ,$ for any $A $, $\overline{A}\in S\left(
d\right) $ with $A\geq \overline{A}.$
\end{enumerate}
Then the sequence $\left\{
\frac{S_n}{\sqrt{n}}+\frac{T_n}n\right\} _{n=1}^\infty $, where
$S_n=X_1+\cdot \cdot \cdot +X_n,$ $T_n=Y_1+\cdot \cdot \cdot
+Y_n,$ converges in law to $\xi +\zeta $:
\[
\lim\limits_{n\rightarrow \infty }\hat E\left[ \varphi \left( \frac{S_n%
}{\sqrt{n}}+\frac{T_n}n\right) \right] =\hat E\left[ \varphi
\left( \xi +\zeta \right) \right] ,\quad \forall \varphi \in
C_{b,Lip}\left( R^d\right) ,
\]
where the sublinear functional $G:R^d\times S(d)\rightarrow R$ is
defined by
\[
G\left( p,A\right) :=\hat E\left[ \left\langle p,\zeta
\right\rangle +\frac 12\left\langle A\xi ,\xi \right\rangle
\right] .
\]}

\bibliographystyle{elsarticle-harv}
\bibliography{<your-bib-database>}





\end{document}